\documentclass{amsart}
\usepackage{a4wide,amsfonts,amssymb,amsmath}
\def\proof{{\bf Proof\quad}}
\def\beginpf{\proof}
\def\qed{\hfill\rule{2.2mm}{2.2mm}\vspace{1ex}}
\def\endpf{\qed}

\newtheorem{theorem}{Theorem}[section]

\newtheorem{definition}[theorem]{Definition}
\newtheorem{example}[theorem]{Example}
\newtheorem{lemma}[theorem]{Lemma}
\newtheorem{proposition}[theorem]{Proposition}
\newtheorem{remark}[theorem]{Remark}

\def\w{\omega}

\def\CC{\mathbb C}

\def\NN{\mathbb N}
\def\RR{\mathbb R}
\def\FF{\mathcal F}

\def\LL{\mathcal L}

\def\cD{\mathcal D}

\def\cT{\mathcal T}
\def\cC{\mathcal C}

\newcommand{\HH}{\mathcal H}

\renewcommand{\AA}{\mathcal{A}}

\newcommand{\re}{\mathop{\rm Re}\nolimits}
\renewcommand{\Re}{\re}

\newcommand{\diag}{\mathop{\rm diag}\nolimits}

\newcommand{\Tt}{(T(t))_{t \ge 0}}

\newcommand{\supp}{\mathop{\rm supp}\nolimits}

\def\text{\mbox}

\begin{document}

\title{Zero-class admissibility of observation operators}
\author{Birgit Jacob}
\address{Institut f\"ur Mathematik,
Universit\"at Paderborn, Warburger Stra\ss e 100, 33098 Paderborn, Germany}
\email{jacob@math.uni-paderborn.de}
\author{Jonathan R.~Partington}
\address{School of Mathematics,
University of Leeds,
Leeds LS2 9JT, U.K.}
\email{J.R.Partington@leeds.ac.uk}
\author{Sandra Pott}
\address{Institut f\"ur Mathematik,
Universit\"at Paderborn, Warburger Stra\ss e 100, 33098 Paderborn, Germany and Department of Mathematics, University of Glasgow,
Glasgow, G12 8QW, U.K. }
\email{s.pott@maths.gla.ac.uk}

\begin{abstract}
An admissible observation operator is zero-class admissible if the norm of the output map tends to zero as the time tends to zero. Sufficient and necessary conditions for zero-class admissibility of observation operators are developed and a modified Weiss condition is studied. It is shown that the
modified Weiss condition is in general necessary, but not sufficient for zero-class admissibility. For several important classes of $C_0$-semigroups it is proved that the modified Weiss condition is indeed equivalent to zero-class admissibility. The methods are illustrated by certain PDE examples.
\end{abstract}


\maketitle

\section{Introduction}\label{introduction}

Consider the Cauchy problem
\begin{eqnarray}\label{system}
\dot x(t) &=& Ax(t), \qquad t \ge 0,\nonumber\\
y(t) &=& Cx(t).
\end{eqnarray}
where the state $x(t)$ lies in a Hilbert space $H$  and the output $y(t)$ lies
in a Hilbert space $Y$ for each time $t \ge 0$. Here $A$ and $C$ are linear operators that may be unbounded,
but $A$ is the generator of a $C_0$-semigroup $(T(t))_{t \ge 0}$ on $H$ and $C$ maps $D(A)$, the domain of $A$,
boundedly into $Y$. In order to guarantee that the output function lies locally in $L^2$ we impose the following
condition.

\begin{definition}
\label{def:itaoo}
The operator $C \in \LL(D(A),Y)$ is a {\em finite-time admissible observation operator for $\Tt$},
if for each
$\eta >0$
there is a constant $K_\eta>0$ such that
\begin{equation}\label{eqn1}
\int_0^\eta \|CT(t)x\|^2 \, dt \le K^2_\eta \|x\|^2, \qquad x \in D(A).
\end{equation}
Likewise, $C$ is an {\em infinite-time admissible observation operator}, if in addition the constants $K_\eta$ are uniformly bounded.
\end{definition}

For exponentially stable systems, the two notions are equivalent. The notion of admissible observation
operators is well studied in the literature, see for example \cite{jp}, \cite{staffans}, \cite{tuwe}, \cite{we89}.
In applications, quite often the observation operator belongs to the zero class
of admissible observation operator, see Section \ref{sectionadm} for an example.

\begin{definition}
The finite-time admissible observation operator $C$ is said to belong to the {\em zero class
of admissible observation operators for $\Tt$} ($C$ is  zero-class admissible), if the best constant $K_\eta$, given by (\ref{eqn1}), satisfies
$K_\eta \to 0$ as $\eta \to 0$.
Likewise, $C$ is  {\em infinite-time zero-class admissible (ITZCA)}, if in addition the constants $K_\eta$ are uniformly bounded.

\end{definition}

The class of zero-class admissible observation operators was first introduced in \cite{XLY}, in order to
provide conditions for exact observability of semigroup systems. One of the main results in \cite{XLY} shows that if the observation operator is
zero-class admissible and
 the system is exact observable, then the semigroup generated by $A$
is left invertible. If additionally  the residual spectrum of $A$ is empty, then $A$ generates a $C_0$-group.
Clearly, bounded observation operators belong to the zero class of admissible observation operators. In Section \ref{sectionadm} an example of an admissible but not zero-class admissible observation operator is given. Further, we show that if $A$ generates
an analytic semigroup, $S\in\mathcal{L}(H,Y)$ and $C=S(-A)^\alpha$, for some $\alpha\in(0,\frac{1}{2})$, then $C$ is an admissible observation operator for $\Tt$ of zero class. For positive definite operator $A$ with compact resolvent this result can be found in \cite{XLY}.

A necessary condition for admissibility is the {\em Weiss condition} (A1): There exists a constant $m>0$ such that
\[
 \|C(sI-A)^{-1}x\| \le \frac{m\|x\|}{\sqrt{\re s}}, \qquad x \in H, \quad s \in \CC_+.
\]
The Weiss condition is actually equivalent to infinite-time admissibility for several classes of systems.
However, it has been shown that in general the Weiss condition is not equivalent to infinite-time admissibility.
We refer the reader to the survey article \cite{jp} for more information on the Weiss condition.
In this article we introduce a modified Weiss condition in order to characterize zero-class admissibility.

We proceed as follows. In Section \ref{section2} we review some known results on admissibility. The main results of this paper are contained
in Section \ref{sectionadm}. We give necessary and sufficient condition for zero-class admissible observation operator.
In particular a modified Weiss condition is introduced and studied. The obtained results are illustrated by certain PDE examples. Finally, in Section \ref{section4} we give two examples showing that the
modified Weiss condition is not sufficient for zero-class admissibility.


\section{Equivalent conditions for admissibility}\label{section2}

In this section we review some known results on admissibility and we pay particular care to the constants involved.
In \cite{JPP03}, it was shown that for bounded $C_0$-semigroups, the following conditions (among others)
are equivalent:
\begin{itemize}
\item (A1) There exists a constant $m>0$ such that
\[
\|C(sI-A)^{-1}x\| \le \frac{m\|x\|}{\sqrt{\re s}}, \qquad x \in H, \quad s \in \CC_+;
\]
\item (A2a) There exists a constant $K>0$ such that
\[
\| \tau^{-1/2} \int_0^\tau e^{i\omega t} CT(t) x \, dt \| \le K \|x\|, \qquad x \in D(A), \quad \tau>0, \quad \omega \in \RR;
\]
\item (A2b) There exists a constant $K>0$ such that
\[
\| \tau^{-1/2} \int_\tau^{2\tau} e^{i\omega t} CT(t) x \, dt \| \le K \|x\|, \qquad x \in D(A), \quad \tau>0, \quad \omega \in \RR.
\]
\end{itemize}
Moreover, these conditions are all implied by admissibility and under certain circumstances (e.g.~when $\Tt$ is a contraction semigroup and $Y$ is finite-dimensional, or when $\Tt$ is an exponentially stable right-invertible semigroup, or when $\Tt$ is analytic and $(-A)^{1/2}$ is admissible, see \cite{japa}, \cite{leme} and \cite{we}.), they
are equivalent conditions to admissibility. Note that condition (A1) is the Weiss condition introduced in Section \ref{introduction}.

More precisely, we have the following result. We recall that $S$ is an {\em $\omega$-sectorial operator if} $S$ is a closed linear operator $S$ in a Hilbert space $H$ satisfying $\sigma(S)\subset
\{s\in\mathbb C\backslash\{0\}\mid |\mbox{arg}\,\lambda|<\omega\}$
and for every $\nu\in(\omega,\pi)$ we have
\[ \sup\{ \|\lambda(\lambda I-S)^{-1}\|\, \mid \,|\mbox{arg}\,\lambda|\ge \nu\} <\infty.\]
An operator $A$ generates a bounded analytic semigroup if and only if $-A$ is a densely defined $\omega$-sectorial operator with $\omega\in(0,\frac{\pi}{2})$.

\begin{proposition}\label{sec:equiv}
Let $N\in \mathbb N$, $M\ge 1$, $\alpha>0$ and $\omega\in(0,\frac{\pi}{2})$, let $\Tt$ be a bounded $C_0$-semigroup with infinitesimal generator $A$ and  $C\in \mathcal{L}(D(A),Y)$. We assume further that one of the following conditions holds:
\begin{enumerate}
\item $\Tt$ is a contraction semigroup and dim$\,Y\le N$.
\item $-A$ is an $\omega$-sectorial operator and $(-A)^{1/2}$ is infinite-time admissible.
\item $\Tt$ is an exponentially stable right-invertible semigroup with $\|T(t)\|\le M e^{-\alpha t}$.
\end{enumerate}
Then $C$ is an infinite-time admissible observation operator if and only if Property (A1) holds.

Denoting by $K:=\sup_{\eta>0} K_\eta$ the best constant in Definition \ref{def:itaoo} and by $m$ the best constant in Property (A1), there exists a constant $c>0$, only dependent on $N$, $M$, $\alpha$ and $\omega$, such that $K\le c m$.
\end{proposition}

\beginpf
It is well-known that  $C$ is an infinite-time admissible observation operator if and only if Property (A1) holds, see \cite{japa}, \cite{leme} and \cite{we}.
Suppose that $K$ is not bounded by an absolute multiple of $m$. Then there exists for each $k\in\mathbb N$,  a semigroup $(T_k(t))_{t \ge 0}$ on $H_k$ with infinitesimal generator $A_k$, and an observation operator
$C_k: D(A_k) \to Y$, such that  $(T_k(t))_{t \ge 0}$ and $C_k$ satisfy the assumption of the proposition, $m_k=1$ and  $(K_k)_{k \ge 1}$ is unbounded. Here $K_k:=\sup_{\eta>0} (K_k)_\eta$ denotes the best constant of $C_k$ in Definition \ref{def:itaoo} and  $m_k$ is the best constant in Property (A1).
Without loss of generality we may assume that each semigroup $(T_k(t))_{t \ge 0}$ satisfies the same Condition 1, 2, 3 or 4.

Then
we form the product semigroup $\Tt$ on the $\ell^2$ direct sum $H$ of the spaces $H_k$, with generator $A=\diag(A_1,A_2,\ldots)$, and
for every $\ell^1$ sequence $(\alpha_k)$ we consider the observation operator
$C$ defined by
\[
C(x_1,x_2,\ldots)=\sum_{k=1}^\infty \alpha_k C_k x_k.
\]
Clearly, we may choose $(\alpha_k)$ such that $\sum_{k=1}^\infty |\alpha_k|=1$ but $(K_k \alpha_k)$ is unbounded.

Since each semigroup $(T_k(t))_{t \ge 0}$ satisfies the same Condition 1, 2 or 3, the product
semigroup satisfies the same condition. Thus for the product semigroup we have the equivalence of
infinite-time admissibility and Property (A1).

Now $C(sI-A)^{-1}x= \sum_{k=1}^\infty \alpha_k C_k (sI-A_k)^{-1}x_k$, so that we have
\[
\|C(sI-A)^{-1}x\| \le  \sup_{k} \|C_k(sI-A_k)^{-1}x_k\| \le \frac{\|x\|}{\sqrt{\re s}},
\]
and $CT(t)x=\sum_{k=1}^\infty \alpha_k C_k T_k(t)x_k$, so for each $k$
\[
\sup_{\|x\| \le 1} \|CT(\cdot)x\|_{L^2(0,\infty)} \ge   \sup_{\|x_k\| \le 1}   |\alpha_k|\|C_k T_k(\cdot)x_k\|_{L^2(0,\infty)} = |\alpha_k| K_k.
\]
Thus
$C$ satisfies the resolvent condition but not the admissibility condition.
This is a contradiction to our assumption and thus the proposition is proved.
\endpf


\section{Equivalent conditions for zero-class admissibility}\label{sectionadm}

In this section we develop equivalent conditions for zero-class admissibility.
Consider the following conditions, which we shall see are equivalent to infinite-time zero-class admissibility in some circumstances.
Note that condition (B1) is a modification of the Weiss condition (A1).

\begin{itemize}
\item (B1) For each
 $r>0$  there exists a constant $m_r>0$ such that the   $m_r$ are uniformly bounded,
$m_r \to 0$   as $r \to \infty$, and
\[
\|C(sI-A)^{-1}x\| \le \frac{m_{\re s}\|x\|}{\sqrt{\re s}}, \qquad x \in H, \quad s \in \CC_+.
\]
Clearly we can, and will, assume without loss of generality that $m_r$ is a decreasing function of $r$.
\item (B2a) For each $\tau>0$ there exists a constant $K_\tau>0$ such that the   $K_\tau$
are uniformly bounded, $K_\tau \to 0$ as $\tau \to 0$, and
\[
\| \tau^{-1/2} \int_0^\tau e^{i\omega t} CT(t) x \, dt \| \le K_\tau \|x\|, \qquad x \in D(A), \quad \omega \in \RR;
\]
\item (B2b) For each $\tau>0$ there exists a constant $K_\tau>0$ such that the   $K_\tau$
are uniformly bounded, $K_\tau \to 0$ as $\tau \to 0$, and
\[
\| \tau^{-1/2} \int_\tau^{2\tau} e^{i\omega t} CT(t) x \, dt \| \le K_\tau \|x\|, \qquad x \in D(A), \quad \omega \in \RR.
\]
\end{itemize}
Again we can, and will, assume without loss of generality that $K_\tau$ is a increasing function of $r$.

It is easily seen that (B2a) and (B2b) are equivalent, first by writing $\int_0^\tau$ as $\sum_{n=0}^\infty \int_{2^{-n-1}\tau}^{2^{-n}\tau}$, and second
by writing $\int_\tau^{2\tau}=\int_0^{2\tau}-\int_0^\tau$, in each case making the obvious estimates as done for the equivalence between (A2a) and (A2b) in
\cite{JPP03}. We may therefore write condition (B2) to refer to either of these conditions.\\

\begin{remark}\label{rem:itzcab2}
Note that, by the Cauchy--Schwarz inequality, infinite-time zero-class admissibility implies condition (B2) immediately.
\end{remark}

\begin{theorem}\label{thm:b1b2}
Let $\Tt$ be a bounded $C_0$-semigroup with infinitesimal generator $A$ and let $C \in \LL(D(A),Y)$ be infinite-time
admissible. Then conditions (B1) and (B2) are equivalent.
\end{theorem}
\beginpf
(B2) $\implies$ (B1). For $s=\sigma+i\omega \in \CC_+$ and $x \in D(A)$ we have
\[
C(sI-A)^{-1}x = \int_{y=0}^\infty \sigma e^{-\sigma y} \int_{t=0}^y CT(t)xe^{-i\omega t} \, dt \, dy,
\]
as  in \cite[p.~321]{JPP03}, and so
\[
\|C(sI-A)^{-1}x\| \le \|x\| \int_0^\infty K_y \sigma e^{-\sigma y} y^{1/2} \, dy,
\]
where $K_y$ is defined in condition (B2a).
Let $y=v/\sigma$, so that
\[
\int_0^\infty K_y \sigma e^{-\sigma y} y^{1/2} \, dy = \sigma^{-1/2} \int_0^{\sqrt\sigma} v^{1/2}e^{-v}K_{v/\sigma} \, dv
+ \sigma^{-1/2} \int_{\sqrt\sigma}^\infty v^{1/2}e^{-v}K_{v/\sigma} \, dv
\]
which is clearly bounded by a constant times $\sigma^{-1/2}$ and
is $o(\sigma^{-1/2})$ as $\sigma \to \infty$, given that $K_\eta$ is a bounded function of $\eta$, which tends to 0 as $\eta\to 0$.\\

(B1) $\implies$ (B2).
We adapt an argument from \cite[pp.~319--320]{JPP03}, paying particular care to the constants involved.
Note first that, with $J(s)=C(sI-A)^{-1}$ we have
\[
\int_0^\tau CT(t)x e^{-(\rho/\tau+i\omega)t}\, dt = J(\rho/\tau+i\omega)x-e^{-\rho}e^{-i\omega \tau}J(\rho/\tau+i\omega)T(\tau)x,
\]
and so
\begin{equation} \label{eq:eqrhotau}
\left\| \int_0^\tau CT(t)x e^{-(\rho/\tau+i\omega)t} \, dt \right\|
\le m_{\rho/\tau}(1+e^{-\rho}K)\frac{\tau^{1/2}}{\rho^{1/2}}\|x\|,
\end{equation}
where $K$ is the norm bound of the semigroup $\Tt$.\\

Let $0 < \rho <\frac12$. We  have from \cite[Lem.~2.4]{JPP03}
that $\chi_{[0,\tau)}=\sum_{I \in \cD_\tau} a_I \psi_I$, where $\cD_\tau$ is the collection of half-open dyadic subintervals
$I \subseteq [0,\tau)$, with left endpoint $l(I)$, and $a_I \ge 0$ for all $I$,
with $\psi_I(t)=|I|^{-1/2} \chi_I e^{-\rho|I|^{-1}(t-l(I))}$, the sum converging uniformly on $[0,\tau)$.
Moreover, $\sum_{I \in \cD_\tau} a_I \le M\tau^{1/2}$, where the constant $M$ can be taken
independent of $\tau$.\\

Using this and (\ref{eq:eqrhotau}), we have
\begin{eqnarray*}
\left\| \int_0^\tau e^{i\omega t}CT(t)x \, dt \right\| &\le& \sum_{I \in \cD_\tau} a_I \left\| \int_I \psi_I(t) e^{i\omega t} CT(t) x \, dt \right\| \\
& \le & \sum_{I \in \cD_\tau} a_I |I|^{-1/2} \left \| \int_0^{|I|} e^{-\rho|I|^{-1} t} e^{i\omega (t+l(I))} CT(t)T (l(I))x \, dt \right \| \\
&\le& K (1+e^{-\rho}K)\rho^{-1/2} \|x\| \sum_{I \in \cD_\tau} a_I m_{\rho/|I|}\\
&\le& K  (1+e^{-\rho}K)\rho^{-1/2}\|x\|m_{\rho/\tau} M\tau^{1/2},
\end{eqnarray*}
and the result follows since $m_{\rho/\tau} \to 0$ as $\tau \to 0$.
\endpf

We are now in a position to present the main result of this section.

\begin{theorem}
\label{thm:heregoes}
Let $\Tt$ be a bounded $C_0$-semigroup with infinitesimal generator $A$ and let $C \in \LL(D(A),Y)$.
Assume that one of the following conditions holds:
\begin{itemize}
\item $\Tt$ is a contraction semigroup and $Y$ is finite-dimensional.
\item $\Tt$ is an exponentially stable right-invertible semigroup.
\item $\Tt$ is an analytic semigroup and $(-A)^{1/2}$ is infinite-time admissible.
\end{itemize}
Then $C$ is an  infinite-time zero-class admissible observation operator (ITZCA)
if and only if condition
(B1) holds.
\end{theorem}
\beginpf
(B1) $\implies$ (ITZCA).
Take $\eta>0$ and write $\lambda=1/\eta$. Then (B1) implies that
\[
\|C((s+\lambda)I-A)^{-1}x\| \le \frac{m_{\re s+\lambda}\|x\|}{\sqrt{\re s+\lambda}} \le
\frac{m_{\lambda}\|x\|}{\sqrt{\re s}},\qquad s \in \CC_+, \quad x \in H.
\]
This implies infinite-time admissibility for $C$ with respect to the semigroup $(e^{-\lambda t}T(t))_{t \ge 0}$
generated by $A-\lambda I$. Using (B1) we have $\lim_{\lambda \rightarrow\infty}m_{\lambda} = 0$.
Then, for $x \in D(A)=D(A-\lambda I)$,
\[
\int_0^\eta \|CT(t)x\|^2 \, dt \le e^2\int_0^\eta \|Ce^{-\lambda t}T(t)x\|^2 \, dt \le e^2 \int_0^\infty \|Ce^{-\lambda t}T(t)x\|^2 \, dt
\le e^2 C_\lambda^2 \|x\|^2.
\]
By Proposition \ref{sec:equiv} the constant $C_\lambda$ is bounded by an absolute multiple of $m_\lambda$. Thus $C_\lambda \to 0$ as $\lambda \to \infty$. Recalling that $\eta=1/\lambda$,
we have the result.

(ITZCA) $\implies$ (B1). This follows from Remark \ref{rem:itzcab2} and Theorem \ref{thm:b1b2}.
\endpf

\begin{remark}
Suppose that $(T(t))_{t \ge 0}$ has a Riesz basis $(\phi_n)_{n \ge 1}$ of eigenvectors, with eigenvalues $(\lambda_n)_{n \ge 1}$.
For an observation operator $C: D(A) \to \CC$ satisfying $C\phi_n=c_n$ ($n \ge 1$), it is easily verified that
$C$ satisfies the  condition (B1), which is here equivalent to (ITZCA), if and only if
\[
\int_{\CC_+} \frac{d\mu(\lambda)}{|s+\lambda|^2} \le \frac{m_{\re s}}{\re s}, \qquad \hbox{for} \quad s \in \CC_+,
\]
where $\mu=\sum_{n=1}^\infty |c_n|^2 \delta_{-\lambda_n}$
and the constants $m_r>0$ are uniformly bounded with $m_r \to 0$ as $r \to \infty$.
Here $\delta_{-\lambda_n}$ denotes a Dirac measure at $-\lambda_n$.

By standard estimates, this property is seen to be equivalent to saying that for a Carleson square
\[
Q_{r,\omega}=\{x+iy \in \CC:\, x \in [0,r], \, y \in [\omega-r/2 , \omega+r/2] \}
\]
 with $r>0$ and $\omega \in \RR$,
one has that $\mu(Q_{r,\omega})/r \le M_r$, where $M_r>0$ is uniformly bounded and $M_r \to 0$ as $r \to \infty$ (cf. \cite{horu,we88}).

Easy examples show:

(i)~This property is strictly stronger than the property of being a {\em Carleson measure\/} (i.e.,
$\mu(Q_{r,\omega})/r$ is bounded;  equivalently,
the canonical embedding $J$ from $H^2(\CC_+)$ to $L^2(\CC_+,\mu)$ is bounded):
e.g., take $\sum_{n=1}^\infty \delta_n$.

(ii)~The property is strictly weaker than the property of being a {\em vanishing Carleson measure\/} (i.e.,  for all $\varepsilon>0$
there is a compact subset $K \subset \CC_+$ such that
$\mu(Q_{r,\omega})/r < \varepsilon$ whenever $\frac{r}{2}+i\omega \not\in K$; equivalently,
$J$ is a compact operator, cf. \cite{garnett}):
e.g., take $\sum_{n=1}^\infty n^{-2}\delta_{1/n}$.
\end{remark}

Next we give some sufficient conditions for zero-class admissibility.

\begin{theorem}
Let $(T(t))_{t\ge 0}$ be an analytic and bounded $C_0$-semigroup with infinitesimal generator $A$ and let $C=S(-A)^\alpha$, where $S$ is a
linear bounded operator from $H$ to $Y$ and $\alpha\in(0,\frac{1}{2})$. Then $C$ is a zero-class admissible observation operator.
\end{theorem}
\beginpf
Using \cite[Theorem 6.13, p.~74]{pazy} there exists a constant $M>0$ such that
\[ \int_0^{\eta} \|CT(t)x\|^2dt \le \|S\|^2 \int_0^{\eta} M^2 t^{-2\alpha}\,dt \|x\|^2 \le  K_\eta \|x\|^2, \qquad \hbox{for} \quad x \in D(A),
\]
where $K_\eta \to 0$ as $\eta \to 0$.
\endpf

A further sufficient condition for zero-class admissibility, which is an analogue of \cite[Prop.~2.1]{XLY} for
normal subgroups (and hence also those with a Riesz basis of eigenvectors) is the following, which includes
also the case when $-A$ is $\omega$-sectorial (where we can take $\beta=1$).

\begin{theorem}
Let $(T(t))_{t \ge 0}$ be a normal semigroup with generator $A$ such that $\sigma(-A)$ is contained
in a region $\{z=v+i w \in \CC: \, v > 0, \, |w| \le a+bv^\beta\}$ for some $a$, $b>0$ and $\beta > 0$.
If $C=S(-A)^\alpha$ for some bounded linear operator $S: H \to Y$ and $\alpha \in (0,\min(\frac12,\frac{1}{2\beta}))$, then $C$ is zero-class admissible.
\end{theorem}

\beginpf
By standard functional calculus arguments for normal operators, we may conclude that, for
$t >0$ and $x \in D(A)$,
\[
\|CT(t)x\| \le \|S\| \sup_{z \in \sigma(-A)} |z|^\alpha |e^{-zt}| \|x\|.
\]
This is bounded by $c_1\max_{v \ge 0}\max(e^{-vt},v^\alpha e^{-vt},v^{\alpha\beta}e^{-vt})\|S\|\,\|x\|$ for some $c_1>0$,
which in turn is bounded by $c_2 \max (1,t^{-\alpha},t^{-\alpha\beta})\|S\|\,\|x\|$ for some $c_2>0$: both $c_1$ and $c_2$
can be taken to be independent of $t$.
Thus, if $\alpha \in (0,\min(\frac12,\frac{1}{2\beta}))$, we see that
\[
\int_0^\eta \|CT(t)x\|^2 \le K_\eta \|x\|^2, \qquad \hbox{for} \quad x \in D(A),
\]
where $K_\eta \to 0$ as $\eta \to 0$.
\endpf

In \cite{zwart}, Zwart proved a sufficient condition for infinite-time admissibility. Actually his proof shows zero-class admissibility.

\begin{theorem}\label{zwart}
Let $(T(t))_{t\ge 0}$ be a bounded $C_0$-semigroup with infinitesimal generator $A$ and let $C\in \mathcal{L}(D(A),Y)$. If there exists a constant
$K>0$ and a monotonically increasing function $g:(0,\infty)\rightarrow  (0,\infty)$ satisfying
\[ \sum_{n=-\infty}^\infty g(\alpha^n)^{-2}<\infty\quad\mbox{ for some }\alpha >1,\]
such that
\begin{equation}\label{hans} \|C(sI-A)^{-1}\| \le \frac{m}{g({\rm Re}\,s)\sqrt{{\rm Re}\,s}},\quad s\in \mathbb C_+,
\end{equation}
then $C$ is a zero-class infinite-time admissible observation operator.
\end{theorem}

A suitable choice for $g$ is $g(t)=(\log(2+t))^\alpha$ for $\alpha>1/2$. Note that the condition (\ref{hans}) is stronger than the modified Weiss condition ($m_r=m/g({\rm Re}\,s)$). We conclude this section by some examples.

\begin{example}
\label{heat1} {\rm We study the one dimensional heat equation on the interval
$[0,1]$ with Neumann boundary conditions and Dirichlet boundary observation, which is given by
\begin{eqnarray*}
\frac{\partial z}{\partial t}(\xi,t) &=& \frac{\partial^2 z}{\partial \xi^2}(\xi,t),\qquad
           \xi\in(0,1), t\ge 0,\\
\frac{\partial z}{\partial\xi}(0,t) &=& 0, \quad  \frac{\partial z}{\partial\xi}(1,t)=0,\quad t\ge 0,\\
z(\xi,0) &=& z_0(\xi),\qquad \xi\in(0,1),\\
y(t) &=& z(0,t),\quad t\ge 0.
\end{eqnarray*}
This p.d.e.~can be written equivalently in the form (\ref{system}) with $H=L^2(0,1)$,
$A$ is given by
\[ A \phi_n = \lambda_n \phi_n,\quad n\in\mathbb N_0,\]
with $\phi_n(x):=\sqrt{2}\cos(n\pi x)$ and $\lambda_n:=-\pi^2n^2$, and $C$ is defined by
$C\phi_n =\sqrt{2}$. Note that $(\phi_n)_{n\in\mathbb N_0}$ is an orthonormal basis of
$H$. An easy calculation shows that $C$ is a zero-class infinite-time admissible observation operator.}
\end{example}

\begin{example}
{\rm The one dimensional undamped wave equation with Dirichlet boundary conditions and Neumann observation
can be described by
\begin{eqnarray*}
\frac{\partial^2 z}{\partial t^2}(\xi,t) &=& \frac{\partial^2 z}{\partial \xi^2}(\xi,t),\qquad
           \xi\in(0,1), t\ge 0,\\
z(0,t) &=& 0, \quad   z(1,t)=0,\quad t\ge 0,\\
z(\xi,0) &=& z_0(\xi),\quad \frac{\partial z}{\partial t}(\xi,0) = z_1(\xi),\quad \xi\in(0,1),\\
y(t) &=& \frac{\partial z}{\partial \xi}(0,t),\quad t\ge 0.
\end{eqnarray*}
The partial differential equation can be written in form (\ref{system}) by choosing $H=D(A_0^{1/2})\times L^2(0,1)$,
$A=\left(\begin{smallmatrix} 0 & I\\ -A_0 &0 \end{smallmatrix}\right)$ with $D(A)=D(A_0)\times D(A_0^{1/2})$ and
$C w= (1 \,\,\, 0)\frac{\partial w}{\partial \xi}(0)$. Here $A_0h= \frac{d^2 h}{d \xi^2}$ with
$D(A_0)=\{ h\in H^2(0,1)\mid h(0)=h(1)=0\}$. An easy calculation shows that $A$ has the eigenvalues $\lambda_n=in\pi$, $n\in\mathbb Z$,
$n\not=0$, and the corresponding normalized eigenvectors are given by $\phi_n(x)=\lambda_n^{-1}\left(\begin{smallmatrix} \sin (n\pi x)\\
\lambda_n \sin (n\pi x) \end{smallmatrix}\right)$; moreover $(\phi_n)_{n\in\mathbb Z, n\not=0}$ is an orthonormal basis of
$H$. It is easy to see that $C$ is an admissible observation operator, but not a zero-class admissible observation operator.
}
\end{example}

\section{Counterexamples}\label{section4}

In this section we give two examples showing that the modified Weiss condition is in general not sufficient for zero-class admissibility. Note that these examples even show that the modified Weiss condition in general does not imply admissibility.
We start with the case of an analytic, exponentially stable semigroup and scalar-valued outputs.

\begin{theorem}\label{thm:4.1} There exists an analytic, exponentially stable semigroup $(\cT(t))_{t \ge 0}$  with infinitesimal generator $\AA$ on a
separable Hilbert
space $\HH$ and an observation operator $\cC: \cD(\AA) \rightarrow \CC$ such that
\begin{enumerate}
\item Condition (B1) holds;
\item there exist a sequence of positive numbers $(c_N)$ and a sequence $(x_N)$ in $\HH$ such that
$c_N \to \infty$, $\|x_N\|=1$ and
$$
  \int_0^1 \|\cC \cT(t) x_N \|^2 \ge c_N \|x_N\|^2   \text{ for all } N \in \NN.
$$
\end{enumerate}
That means, $(\AA,\cC)$ satisfies (B1), but
is not finite-time admissible for any time $\eta>0$, and in particular not zero-class admissible.
\end{theorem}

The proof relies on the construction in \cite{JZ}, but some care has to be taken to adapt this to the setting of
zero-class admissibility. Let $H$ be an
infinite-dimensional separable Hilbert space, let $(\phi_n)$ be a non-Besselian conditional basis of $H$
with $\inf_{n \in \NN} \|\phi_n\| >0$, let $\mu_n = - 4^n$ for $n \in \NN$, and let $A$ be densely defined on
$H$ defined by
$A \phi_n = \mu_n \phi_n$ for $n \in \NN$. As shown \cite{JZ}, $A$ is the infinitesimal generator of
an analytic, exponentially stable semigroup $\Tt$. For $x = \sum_{n=1}^\infty \alpha_n \phi_n
\in D(A)$ and $N \in \NN$,  we define
$$
   C_{N} x = \sum_{n=1}^N \sqrt{-\mu_n} \alpha_n \phi_n.
$$
As in Proposition 3.1 of \cite{JZ}, one shows that the $C_{N}$ are uniformly $A$-bounded for all $N$.
Here is our main estimate.
\begin{lemma}\label{lem:4.2}
There exists a sequence of positive numbers $(c_{N})_{N \in \mathbb N}$ with
$c_{N} \to \infty$ for $N \to \infty$, a sequence of vectors $(x_{N})_{N \in \mathbb N}$, $\|x_{N}\|=1$,
and sequences $(m_{N,r})_{N \in \mathbb N}$, with $m_{N,r} \to 0$ as $r \to 0$ for each
$N \in \NN$ and $m_{N,r}$ uniformly bounded in $N$ and $r$,
  such that
$$
   \|C_{N}(sI-A)^{-1}\| \le m_{N,\re s} \frac{1}{\sqrt{\re s}},\quad s\in\mathbb C_+, N\in\mathbb N,
$$
and
$$
   \int_0^{1} |C_{N} T(t) x_{N}|^2 dt \ge c_{N} \| x_{N}\|^2,\quad N\in\mathbb N.
$$
\end{lemma}
{\bf Proof of Lemma \ref{lem:4.2}\quad}
By the same calculation as in \cite{JZ}, Prop.~3.2,
we have for $x \in H$, $\|x\|=1$, and for each $N \in \NN$:
$$
   \sqrt{\re s}|C_{N} (s-A)^{-1}x|  \le \kappa \sqrt{\re s} \sum_{n=1}^N \frac{2^n}{\re s + 4^n}
   \le 2 \kappa \sqrt{\re s} \sum_{n=1}^{2^N} \frac{1}{\re s + n^2},
$$
where $\kappa$ is an absolute constant only depending on the sequence $(\phi_n)$.

Letting $m_r =2 \kappa r^{1/2} \sum_{n=1}^{\infty} \frac{1}{r + n^2}$ and
$m_{N,r}=2 \kappa r^{1/2} \sum_{n=1}^{2^N} \frac{1}{r + n^2}\le 2\kappa\arctan\frac{2^N}{\sqrt{r}}$, we find that $m_{N,r} \le m_r \le \kappa \pi$ for all $r >0$
as in
\cite{JZ}, $\lim_{r \to \infty} m_{N,r}=0$ for each $N \in \NN$, and
$$
 |C_{N} (sI-A)^{-1}x| \le \frac{m_{N,\re s}}{\sqrt{\re s}}.
$$
Consequently, $(A, C_{N})$ satisfies (B1) for each $N \in \NN$, with a uniform bound
on the $m_{N,r}$.

On the other hand, since $(\phi_n)$ is non-Besselian, we can find a sequence of positive numbers
$(\tilde c_{N})$, $\tilde c_{N} \to \infty$ for $N \to \infty$,
and vectors $x_N = \sum_{n=1}^N \alpha_{N,n} \phi_n$, $\|x_N\|=1$, such that
$$
   \| x_N \|^2 \le \frac{1}{\tilde c_{N}} \sum_{n=1}^N |\alpha_{N,n}|^2.
$$
However, the system $(\sqrt{-\mu_n} e^{\mu_n t})_{n \ge 1}$ is unconditional in $L^2(0,1)$
(see e.~g.~\cite{Nik}, Corollary 4.5.2).
Therefore, there exists a constant $L >0$ independent of $N$ such that
$$
   \int_0^1 |C_{N} T(t) x_{N}|^2 dt
   =  \int_0^1 |\sum_{n=1}^N \alpha_{N,n} \sqrt{-\mu_n} e^{\mu_n t} \phi_n|^2 dt \ge
   L \sum_{n=1}^N |\alpha_{N,n}|^2.
$$
Letting $c_{N} = L \tilde c_{N}$, we find that
$\displaystyle
     \int_0^1 |C_{N} T(t) x_{N}|^2 dt \ge {c_{N}} \|x_{N}\|^2.
$
\endpf

{\bf Proof of Theorem \ref{thm:4.1} \quad}
Let $\alpha\in(0,1/2)$.
We choose a sequence of nonnegative numbers $(\beta_N)$ with
$\sum_{N=1}^\infty |\beta_N|=1$, $(\beta_N^2 c_N^{1-2\alpha})$ unbounded.

Now as in the proof of Proposition \ref{sec:equiv}, form the $\ell^2$-direct sum of countably many copies $H_N$ of the Hilbert space $H$,
$\HH= \bigoplus_{N=1}^\infty H_N$,
let $\AA$ denote the densely defined block-diagonal operator $\diag(A,A, \dots)$ on $\HH$, and let $(\cT(t))_{t \ge 0}$
denote the semigroup generated by $\AA$. Obviously $(\cT(t))_{t \ge 0}$ is still analytic and exponentially stable. Define
$$
   \cC: D(\AA) \rightarrow \CC, \qquad \cC x= \sum_{N=1}^\infty \frac{\beta_N }{c_N^{\alpha}} C_{N} x_N,
$$
for $x= (x_1, x_2, \dots) \in \HH$.
Then for $x= (x_1, x_2, \dots) \in \HH$, and $s \in \CC_+$,
$$
   |\cC(sI-\AA)^{-1} x | = |\sum_{N=1}^\infty  \frac{\beta_N}{c_N^{\alpha}} C_{N} (sI-A)^{-1} x_N |
      \le \sup_{N \in \NN} \frac{1}{c_N^{\alpha}} |C_{N} (sI-A)^{-1} x_N |\le
          \sup_{N \in \NN} \frac{1}{c_N^{\alpha}} \frac{m_{N,\re s}}{\sqrt{\re s}}.
$$
It is easy to see that
$
    M_r = \sup_{N \in \NN} \frac{1}{c_N^{\alpha}} m_{N, r}
$
satisfies $M_r \to 0$ for $r \to \infty$. Thus
$(\AA, \cC)$ satisfies condition (B1).

On the other hand, letting $\tilde x_N = (0, 0, \dots, x_N, 0, 0,\cdots ) \in \HH$,
$$
  \int_0^{1}  |\cC \cT(t) \tilde x_{N}|^2 dt \ge \beta_N^2 c_N^{1-2\alpha} \|\tilde x_{N}\|^2 \text{ for each } N \in \NN
$$
by the lemma. Hence $(\AA, \cC)$ is not finite-time admissible and in
particular not zero-class admissible.
\endpf

Next we study the modified Weiss condition
\begin{equation}\label{hans2} \|C(sI-A)^{-1}\| \le \frac{m}{(\log({\rm Re}\,s+2))^\gamma\sqrt{{\rm Re}\,s}},\quad s\in \mathbb C_+,
\end{equation}
with $\gamma>0$. Theorem \ref{zwart} implies  that condition (\ref{hans2})
is sufficient for zero-class infinite-time admissibility if $\gamma>\frac12$. Next we show that this result
is sharp. More precisely, for $\gamma\in (0,\frac12)$ there exists  an analytic exponentially stable $C_0$-semigroup $(T(t))_{t\ge 0}$  with infinitesimal generator $A$ and an operator $C\in \mathcal{L}(D(A),Y)$
such that $C$ is not finite-time admissible but (\ref{hans2}) holds.

\begin{remark}
We study again the $C_0$-semigroup and the observation operator defined in
the proof of Theorem \ref{thm:4.1} for a particular choice of the non-Besselian basis $(\phi_n)$.
For $0 < \beta < \frac12$ the functions $\phi_k$ given by
$\phi_{2n}(t)=|t|^\beta e^{int}$ and $\phi_{2n+1}(t)=|t|^\beta e^{-int}$,
$t \in (-\pi,\pi)$
($n=0,1,2,\ldots$)
form a non-Besselian bounded basis of $L^2(-\pi,\pi)$ (see
\cite[pp.~351-353]{singer}), and indeed
for $\frac14 < \beta < \frac12$
there is a function
$\sum_{k=0}^\infty \alpha_k \phi_k \in L^2(-\pi,\pi)$ (a multiple of the
function $t \mapsto |t|^{-\beta}$) such that
$\alpha_k$ is the Fourier coefficient  of
the non-$L^2$ function $t \mapsto |t|^{-2\beta}$ with respect to the
orthogonal
system $\{e^{int},e^{-int}, n=0,1,2,\ldots\}$. An easy estimate of
$\int_{-\pi}^\pi |t|^{-2\beta} \exp(\pm int) \, dt$ shows that
$\alpha_k$ is asymptotic to $k^{-1+2\beta}$ and
$\sum_{k=0}^N |\alpha_k|^2$ grows as  $N^{-1+4\beta}$. Thus $c_N$ is asymptotic to
$N^{-1+4\beta}$. This implies that
\[ M_r \le c_1 \sup_{N\in\mathbb N} \frac{\arctan\frac{2^N}{\sqrt{r}}}{N^{(-1+4\beta)\alpha}} \le \frac{c_2}{(\log (r+2))^{(-1+4\beta)\alpha}},\]
for some constants $c_1$, $c_2$ independent on $r$. Note that $\gamma:=(-1+4\beta)\alpha\in(0,\frac12)$.
\end{remark}

Now we give a counterexample for the case of a contraction semigroup, but infinite-dimensional output space.
\begin{theorem}\label{thm:4.3} Let $\HH$ be a separable infinite-dimensional Hilbert space, let
$(\cT(t))_{t \ge 0}$ denote the right shift semigroup on $L^2([0,1]; \HH)$, and let $\AA$ denote the
infinitesimal generator of $(\cT(t))_{t \ge 0}$. Then there exists
an observation operator $\cC: \cD(\AA) \rightarrow \CC$ such that
\begin{enumerate}
\item Condition (B1) holds;
\item there exists a sequence of positive numbers $(c_N)$, $c_N\approx N$, and a sequence $(x_N)$ in $\HH$ such that
$c_N \to \infty$  and
$$
  \int_0^{1} \|\cC \cT(t) x_N \|^2 \ge c_N \|x_N\|^2   \text{ for all } N \in \NN.
$$
\end{enumerate}
That means, $(\AA,\cC)$ satisfies (B1), but
is not finite-time admissible for any time $\eta>0$, and in particular not zero-class admissible.
\end{theorem}

{\em Remark.} This sharpens the counterexample to the Weiss conjecture in \cite{JPP}, where it was shown that the Weiss condition (A1)
in this setting
does not imply infinite-time admissibility.

The proof is based on the counterexample in \cite{JPP}, but again we have to look carefully at the details. Let $H$ be an infinite-dimensional separable Hilbert space, let $e$ be a unit vector in $H$, let $(e_k)$ denote an orthonormal sequence,
and define
$$
    b_N : i \RR \rightarrow H, \quad b_N = \sum_{k=1}^N \psi_k(i \cdot) e_k,
$$
where $(\psi_k)$ is an orthonormal system in $H^2(\CC_+)$ to be specified later.

With the proof of Theorem 2.4 in \cite{JPP}, one sees that the operators given by
\begin{equation}\label{C_N}
     C_N x = \int_{-\infty}^\infty \langle b_N(i \w), x(i\w)\rangle d \w  \quad (x \in \cD(A), h \in H)
\end{equation}
are observation operators, satisfying the estimates
\begin{equation}   \label{eq:b1}
     \|C_N(sI -A)^{-1}\| \le M \|b_N\|_{\mathrm{WBMO}} \frac{1}{\sqrt{\re s}}
\end{equation}
with an absolute constant $M>0$. Here, $\mathrm{WBMO}(H)$ denotes the weak BMO space of $H$-valued functions,
$ \mathrm{WBMO}(H) = \{ b \in L^2(i\RR, H): \sup_{e \in H, \|e\|=1} \|\langle b(\cdot), e \rangle\|_{\mathrm{BMO}} < \infty \}$.

Furthermore, one sees with the calculation in \cite{JPP}, Theorem 2.4, (ii) that for $f \in H^2(\CC_+)$,
\begin{equation} \label{eq:ft}
 C_N T(t) f = \int_0^\infty \FF^{-1}(b_N)(u+t) \FF^{-1}(u)du,
\end{equation}
where $\FF^{-1}$ or $\check{}$ denotes the inverse Fourier transform of functions in
$H^2(\CC_+) \subset L^2(i\RR)$, taken as a
function on $L^2(0, \infty)$.

Let $f_N(i \w) = \frac{1}{N^{1/2}}\sum_{k=1}^N \overline{\psi_k(-i
\w)}$. Then $f_N \in H^2(\CC_+)$, $\|f_N\|_2=1$. Here is the
technical result we require:
\begin{lemma} There exist constants $\delta,K >0$ such that for a suitable choice of $(\psi_k)$,
\begin{enumerate}
\item $\|b_N\|_{WBMO} \le K $ for all $N \in \NN$;
\item $\int_{0}^1 \| \int_0^\infty (\FF^{-1}b_N)(u+t) (\FF^{-1}f)(u)du   \|^2 dt \ge \delta
N$ for all $N \in \NN$.
\end{enumerate}
\end{lemma}
\beginpf Let $\psi: i \RR \rightarrow \CC$, $\psi \in H^2(\CC^+)$, $\|\psi\|_2=1$ such that
\begin{enumerate}
\item $\check \psi$ bounded, $\supp \check{\psi} \subseteq [1,2]$,
\item $|\psi(i \w)|$ is rapidly decreasing,
\item  $\int_{-\infty}^\infty \psi(i\w) =0$,
\item 
$$
\delta = \int_0^1 | \FF^{-1}(|\psi|^2)(-t)|^2 dt =\int_0^1 |\int_0^\infty \check{\psi}(u)
\overline{\check{\psi}(u+t)}du|^2 dt >0
$$
\end{enumerate}
(take for example the analytic part of the Meyer wavelet). Let
$\psi_k = e^{-i \rho_k \cdot} \psi$, $k \in \NN$, where $\rho_k =
2^k$. Then each $\psi_k$ also satisfies properties (1)-(3), and
$(\psi_k)$ is an orthonormal system. Moreover, since the $\check{\psi_k}$ are translates of $\check{\psi}$, supported in
the lacunary family of intervals $(\rho_k +[1,2])$,
all $L^p$ norms, and BMO norm, of linear
combinations of the $\psi_k$ are equivalent (see e.~g.~\cite{St} for
the case of lacunary Fourier series, from which the results fitting
our case are easily deduced). Applying this fact to
scalar-valued functions $i \RR \rightarrow \CC$
 of the form $i \w \mapsto \langle b_N(i \w), e \rangle$, where
$e \in H$, $e= \sum_{k=1}^\infty \alpha_k e_k$ with $\|e\|^2= \sum_{k=1}^\infty \alpha_k e_k=1$, we see that
$$
    \|\langle b_N, e \rangle\|^2_{BMO} = \| \sum_{k=1}^N \alpha_k \psi_k \|^2_{BMO(i \RR)}  \le K\sum_{k=1}^N
    \alpha_k^2 \le K,
$$
where $K$ is an absolute constant only depending on $\psi$ and the
sequence $(\rho_k)$.
Thus the first part of the lemma holds.

For the second part, observe that since $\supp \check{\psi_k}
\subseteq \rho_k +[1,2]$ for each $k$, mixed convolutions
$\check{\psi_k} * \check{\overline{\psi_l(-\cdot)}}$, $l \neq k$, have support
outside $[0,1]$, and
\begin{eqnarray*}
     \int_{0}^1 \| \int_0^\infty \check{b_N}(u+t) \check{f_N}(u)du   \|^2
     dt &= & \frac{1}{N} \int_{0}^1 \| \sum_{k=1}^N \int_0^\infty \check{\psi_k}(u) \overline{\check{\psi_k}(u+t)}du   \|^2
     dt \\
     &= & \frac{1}{N}\int_{0}^1 \| \sum_{k=1}^N \int_0^\infty \check{\psi}(u) \overline{\check{\psi}(u+t)}du   \|^2
     dt = N \delta.
\end{eqnarray*}
This finishes the proof of the lemma.
\endpf

It remains only to prove that the systems $(A, C_N)$ satisfy not only
the Weiss condition (A1), as already proved with (\ref{eq:b1}) and the
lemma, but the stronger condition (B1).
Recall that
$$
      \| C_N(sI -A)^{-1} x \|
      =   \|\int_{-\infty}^{\infty} b_N(i \w) \frac{x(i\w)}{s+ i \w} d
      \w \|.
$$
Since we already have a  uniform
bound over all $N$ in (\ref{eq:b1}), it is sufficient to show that for
each $k$, there
exists a sequence $(m_r)_{r>0}$ such that $m_r \stackrel{r \to
  \infty}{\rightarrow} 0$ and
$$
      |\int_{-\infty}^{\infty} \psi_k(i \w) \frac{x(i\w)}{s+ i \w} d \w|
      \le \frac{m_{\Re s}}{\sqrt{\Re s}} \|x\|.
$$
But this follows from Cauchy-Schwarz and the estimate
\begin{multline*}
   |\FF^{-1}(\psi_k \frac{1}{s+ i \cdot})(t)|= |(\check{\psi_k} * e^{-s
   \cdot})(t)| \le \int_{1+\rho_k}^{2+\rho_k}  \chi_{[0,t)}(u) e^{-{\Re s}(t-u)} |\check{\psi_k}(u)|
   du \\ \le K_1 e^{-{\Re s}t}
   \int_{1+\rho_k}^{\min\{2+\rho_k,t \}}  e^{({\Re s})u} du
   \le K_1\left\{ \begin{matrix} \frac{1}{\Re s} e^{{\Re s}(2+\rho_k
   -t)} & \text{if} & t > 2 + \rho_k, \\
         \frac{1}{\Re s} & \text{if} & 1+ \rho_k \le t < 2 + \rho_k, \\  
         0   & \text{if} & t \le 1+ \rho_k, \\  \end{matrix}   \right.
\end{multline*}
using boundedness of $\check{\psi}$ by a constant $K_1$, which yields
$$
  | \int_{-\infty}^{\infty} \psi_k(i \w) \frac{x(i\w)}{s+ i \w} d \w|
  \le \|x\|_2 \|\FF^{-1}(\psi_k \frac{1}{s+ i \cdot})\|_2 \le K_1
  \|x\|_2 \frac{1}{\Re s} \left(
  1 + \frac{1}{ 2 \Re s}\right)^{1/2}
$$
and therefore our desired result for e.g.~$m_r=2K_1r^{-1/2}$.
We obtain the following lemma for the observation operators $C_N$, defined in (\ref{C_N}):
\begin{lemma}\label{lem:4.4}
 There exists a sequence of positive numbers $(c_{N})_{N \in \NN}$ with
$c_{N} \approx N \to \infty$, a sequence of vectors $(x_{N})_{N \in \NN}$, $\|x_{N}\|=1$,
and sequences $(m_{N,r})_{N \in \NN}$, with $m_{N,r} \to 0$ as $r \to \infty$ for each
$N \in \NN$ and $m_{N,r}$ uniformly bounded in $N,r$,
  such that
$$
   \|C_{N}(sI-A)^{-1}\| \le m_{N,\re s} \frac{1}{\sqrt{\re s}} 
$$
and
$$
   \int_0^{1} |C_{N} T(t) x_{N}|^2 dt \ge c_{N} \| x_{N}\|^2 \quad \text{ for all } N \in \NN.
$$
\end{lemma}

With a weighted diagonal construction similar to the one in Theorem \ref{thm:4.1}, we obtain the desired
counterexample. This finishes the proof of Theorem \ref{thm:4.3}.\\

{\bf Acknowledgements} $\quad$ The authors gratefully acknowledge support from the Royal Society's
{\em International Joint Project\/} scheme. The third author acknowledges support by the German Research Foundation (DFG).

\end{document}